\documentclass[12pt]{amsart}
\usepackage{txfonts}
\usepackage{amssymb}
\usepackage{mathrsfs}
\usepackage{amsmath}
\usepackage{amssymb,latexsym}

\def\Gal{\rm Gal}

\def\Z{\mathbb Z} 

\newtheorem{thm}{Theorem}[section]
\newtheorem{cor}[thm]{Corollary}

\theoremstyle{definition}

\newtheorem{rem}[thm]{Remark}

\numberwithin{equation}{section}

\begin{document}
\baselineskip=17pt

\title
 {$k$th power residue chains of global fields}
\author[S.Hu and Y.Li]{Su Hu and Yan Li}

\address{Department of Mathematical Sciences, Tsinghua University, Beijing 100084,
China}
\email{hus04@mails.tsinghua.edu.cn
\\liyan\_00@mails.tsinghua.edu.cn}

\begin{abstract} In 1974, Vegh proved that if $k$ is a prime and $m$ a positive integer,
there is an $m$ term permutation chain of $k$th power residue for
infinitely many primes [E.Vegh, $k$th power residue chains, J.Number
Theory, 9(1977), 179-181]. In fact, his proof showed that
$1,2,2^2,...,2^{m-1}$ is an $m$ term permutation chain of $k$th
power residue for infinitely many primes. In this paper, we prove
that for any ``possible" $m$ term sequence $r_1,r_2,...,r_m$, there
are infinitely many primes $p$ making it an $m$ term permutation
chain of $k$th power residue modulo $p$, where $k$ is an arbitrary
positive integer [See Theorem ~\ref{main1} below]. From our result,
we see that Vegh's theorem holds for any positive integer $k$, not
only for prime numbers. In fact, we prove our result in more
generality where the integer ring $\Z$ is replaced by any
$S$-integer ring of global fields (i.e. algebraic number fields or
algebraic function fields over finite fields).
\end{abstract}

 \maketitle

\section{Introduction}
Let $K$ be a global field (i.e. algebraic number field or algebraic
function field with a finite constant field). Let $S$ be a finite
set of primes of $K$ (if $K$ is an algebraic number field, $S$
contains all the archimedean primes). Let $A$ be the ring of
$S$-integers of $K$, that is
$$A=\{a\in K|\mathrm{ord}_{P}(a)\geq 0,\forall P\not\in S\}.$$ If $K$ is a number field and $S$ is
the set of the archimedean primes of $K$, then $A$ is just the usual
integer ring $O_K$ of $K$, i.e. the integral closure of $\Z$ in $K$.
 It is well known that $A$ is a Dedekind domain. Let $P$ be a nonzero prime ideal of $A$
and $k$ a positive integer. A sequence of elements in $A$
\begin{equation}\label{1.1}
r_{1},r_{2}\cdots r_{m}
\end{equation}
for which the $\frac{m(m+1)}{2}$ sums
$$\sum_{k=i}^{j}r_{k},1\leq i\leq j\leq m,$$
are distinct $k$th power residues modulo $P$, is called a chain of
$k$th power residue modulo $P$. If
$$r_{i},r_{i+1},\cdots,r_{m},r_{1},r_{2},\cdots,r_{i-1}$$ is a chain
of $k$th power residue modulo $P$ for $1\leq i\leq m$,
then we call (\ref{1.1}) a cyclic chain of $k$th power residue
modulo $P$. If
$$r_{\sigma(1)},r_{\sigma(2)},\cdots,r_{\sigma(m)}$$ is a chain of $k$th power
residues for all permutations $\sigma\in S_m$, then we call
(\ref{1.1}) a permutation chain of $k$th power residue modulo $P$.
These definitions are generalizations of the classical definitions
of $k$th power residue chains of integers modulo a prime number
(see~\cite{vegh}).

Let $k,p$ be prime numbers. In 1974, using Kummer's result on $k$th
power character modulo $p$ with preassigned values, Vegh~\cite{vegh}
proved the following result for $k$th power residue chains of
integers.
\begin{thm}(Vegh~\cite{vegh})
Let $k$ be a prime and $m$ a positive integer. There is an $m$ term
permutation chain of $k$th power residue for infinitely many primes.
\end{thm}
By using the result of Mills (Theorem 3 of~\cite{Mills}), he showed
that this result also holds if the prime $k$ is replaced by other
kinds of integers (for example $k$ odd, $k=4$, or $k=2Q$, where
$Q=4n+3$ is a prime). It should be noted that Gupta~\cite{Gupta}
exhibited quadratic residue chains for $2\leq m\leq 14$ and cyclic
quadratic residues for $3\leq m\leq 6$.

The main result of this paper is the following theorem.
\begin{thm}~\label{main1}Let $k$ and $m$ be arbitrary positive integers.
Let $r_1,r_2,...,r_m$ be a sequence of elements of $A$ such that for
all permutations $\sigma\in S_{m}$,
\begin{equation}the\ m(m+1)/2\ sums
\sum_{k=i}^{j}r_{\sigma(k)}\ (1\leq i\leq j\leq m)\ are\ distinct.
\end{equation}\label{1.2}
Then $r_1,r_2,...,r_m$ is an $m$ term permutation chain of $k$th
power residue for infinitely many prime ideals.
\end{thm}
\begin{rem}By the definition of permutation chain, the condition
(\ref{1.2}) is necessary for $r_1,r_2,...,r_m$ being a permutation
chain of $k$th power residue.
\end{rem}
In Section 2 and 3, we will prove Theorem~\ref{main1} for number
fields and function fields, respectively.
As a corollary, we get the following theorem which is the
generalization of Vegh's Theorem to the case that $k$ is an
arbitrary positive integer and $A$ is any $S$-integer ring of global
fields.
\begin{cor}~\label{main}
  Let $k$ and $m$ be arbitrary positive integers.
In $A$, there is an $m$ term permutation chain of $k$th power
residues for infinitely many prime ideals.
\end{cor}
Proof of Corollary~\ref{main}. Number field case: let $P$ be a prime ideal of
$A$ and $p$ the prime number lying below $P$ and put
\begin{equation}~\label{1} r_{i}=p^{i-1},\ \ i=1,2,\cdots m.
\end{equation}
Function field case: let $t$ be any element of $A$ which is
transcendental over the constant field of $K$ and put
\begin{equation}~\label{2} r_{i}=t^{i-1},\ \ i=1,2,\cdots m.
\end{equation}
It is easy to see $r_1,r_2,...,r_m$ satisfy the
condition of Theorem \ref{main1}.

Our main tool for proving Theorem~\ref{main1} is the following
Chebotarev's density theorem for global fields (Theorem 13.4 of
\cite{Neukrich} and Theorem 9.13A of \cite{Rosen}).
\begin{thm}~\label{Che}\rm{(Chebotarev)} Let $L/K$ be a Galois extension
of global fields with $\Gal(L/K)=H$. Let $C\subset H$ be a conjugacy
class and $S_{K}$ be the set of primes of $K$ which are unramified
in $L$. Then
$$\delta(\{\mathfrak{p}\in S_{K}|(\mathfrak{p},L/K)=C\})=\frac{\#C}{\#H},
$$ where $\delta$ means Dirichlet density. In particular, every
conjugacy class $C$ is of the form $(\mathfrak{p},L/K)$ for
infinitely many places $\mathfrak{p}$ of $K$.
\end{thm}

\section{Proof of the main result for number fields}
Let the set
\begin{equation}~\label{3}
\mathscr{E}=\{\sum_{k=i}^{j}r_{\sigma(k)}|\ \sigma\in S_m,\ 1\leq
i\leq j\leq m \}.
\end{equation}
Define \begin{equation}~\label{4} \mathscr{P}=\{P\ |\ P~\textrm{is a
prime ideal of}~ A\ and\ \exists\
c_{i},c_{j}\in\mathscr{E},c_{i}\not=c_{j}~\textrm{s.t.}~
P|c_{i}-c_{j}\}.\end{equation} It is easy to see that $\mathscr{P}$
is a finite set of prime ideals of $A$ and the elements in
$\mathscr{E}$ modulo $P$ are not equal if $P\not\in\mathscr{P}$.

 Let $\zetaup_{k}$ be a primitive
$k$th roots of unity. Let $L=K(\zetaup_{k},\sqrt[k]{\mathscr{E}})$.
Then $L/K$ is a Kummer extension. By Chebotarev's density theorem,
there are infinitely many prime ideals $P$ in $A$ such that $P$
splits completely in $L$. Let $B$ be the integral closure of $A$ in
$L$ and $\mathfrak{P}$ be a prime ideal of $B$ lying above $P$, then
\begin{equation}~\label{5} \frac{B}{\mathfrak{P}}\cong\frac{A}{P}.
\end{equation} But we have \begin{equation}~\label{6}
c\equiv(\sqrt[k]{c})^{k}~\textrm{mod}~\mathfrak{P},\ \forall
c\in\mathscr{E},
\end{equation} that is $c$ is a $k$th power residue in $B/\mathfrak{P}$. From (\ref{5}), $c$ is also a $k$th power residue in $A/P$.

Let $\mathscr{M}$ be the infinite set of all the prime ideals of $A$
which split completely in $L$. From the above discussions, the
infinite set $\mathscr{M}-\mathscr{P}$ satisfies our requirement.
That is to say all the elements in $\mathscr{E}$ are distinct $k$th
power residues for
 any prime $P$ in $\mathscr{M}-\mathscr{P}$. Hence, $r_1,r_2,...,r_m$ is an $m$ term permutation chain of $k$th
power residue for all the prime ideals $P\in
\mathscr{M}-\mathscr{P}$.

\section{Proof of the main result for function fields}
Let $K$ be a global function field with a constant field
$\mathbb{F}_{q}$, where $q=p^{s}$, $p$ is a prime number.

 1)~If
$(k,p)=1$. We can prove that the sequence $r_1,r_2,...,r_m$ is a
permutation chain of $k$th power residue for infinitely many prime
ideals of $A$ by the same reasoning as in the Section 2.

2)~If $p|k$. Let $k=p^{t}k^{\prime}$ and $(k^{\prime},p)=1$. Let $P$
be a prime ideal of $A$ and $a$ be any element of $A$. Since the
characteristic of the residue field is $p$, it is easy to see that
$a$ is a $k$th power residue modulo $P$ if and only if $a$ is a
$k^{\prime}$th power residue modulo $P$. Since the theorem holds for
$k^{\prime}$ from 1), it also holds for $k$. Thus, we have finished
the proof in this case.

\end{document}